\documentclass{amsart} 

\usepackage{amsmath} 
\usepackage{amssymb} 

\usepackage{amsfonts}

\usepackage[dvipsnames]{xcolor}
\usepackage{tikz}
\usetikzlibrary{cd}
\usepackage{pgfplots}
\pgfplotsset{compat=1.18}
\usepgfplotslibrary{fillbetween}

\begin{document}

\newcommand{\BB}{{\mathbb B}}
\newcommand{\CC}{{\mathbb C}}
\newcommand{\GG}{{\mathbb G}}
\newcommand{\HH}{{\mathbb H}}
\newcommand{\PP}{{\mathbb P}}
\newcommand{\QQ}{{\mathbb Q}}
\newcommand{\RR}{{\mathbb R}}
\newcommand{\TT}{{\mathbb T}}
\newcommand{\ZZ}{{\mathbb Z}}

\newcommand{\cA}{{\mathcal A}}
\newcommand{\cC}{{\mathcal C}}
\newcommand{\Coh}{{\mathcal Coh}}
\newcommand{\cD}{{\mathcal D}}
\newcommand{\cE}{{\mathcal E}}
\newcommand{\cF}{{\mathcal F}}               
\newcommand{\cI}{{\mathcal I}}
\newcommand{\cH}{{\mathcal H}}
\newcommand{\cM}{{\mathcal M}}
\newcommand{\cO}{{\mathcal O}}
\newcommand{\cS}{{\mathcal S}}
\newcommand{\cSets}{{\mathcal Sets}}
\newcommand{\cSub}{{\mathcal Sub}}
\newcommand{\cT}{{\mathcal T}}
\newcommand{\cZ}{{\mathcal Z}}

\newcommand{\bl}{\mbox{bl}} 
\newcommand{\ch}{\mbox{ch}}
\newcommand{\coker}{\mbox{coker}}
\newcommand{\Ext}{\mbox{Ext}}
\newcommand{\Hilb}{\mbox{Hilb}}
\newcommand{\Hom}{\mbox{Hom}} 
\newcommand{\rH}{\mbox{H}} 
\newcommand{\im}{\mbox{Im}}
\newcommand{\length}{\mbox{length}}
\newcommand{\NS}{\mbox{NS}}
\newcommand{\ord}{{\rm ord}}
\newcommand{\rk}{\mbox{rk}}
\newcommand{\re}{{\rm Re}}
\newcommand{\Quot}{{\rm Quot}}
\newcommand{\Sym}{\mbox{Sym}}

\newcommand{\id}{\mbox{id}}
\newcommand{\tr}{\mbox{tr}}
\newcommand{\End}{\mbox{End}}
\newcommand{\Aut}{\mbox{Aut}}
\newcommand{\Perm}{{\rm Perm}}
\newcommand{\sgn}{{\rm sgn}}
\newcommand{\PGL}{{\rm PGL}}
\newcommand{\Pic}{{\rm Pic}}

\newcommand{\Sec}{{\rm Sec}}

\newcommand{\fun}{\rightarrow}
\newcommand{\nt}{\noindent}

\begin{center}

{\bf Two New Extensions of Reider's Theorem on Algebraic Surfaces}

\medskip

Aaron Bertram\footnote{This work was partially supported by a travel grant from the Simons Foundation}, Jonathon Fleck,  Liebo Pan, Joseph Sullivan

\end{center}

\section{Introduction} 

Walls and chambers for the moduli of semi-stable objects  in the derived category of a complex projective manifold have been used to prove
``classical'' theorems about the Brill-Noether theory of curves in surfaces \cite{Ba,BL, BM1, BM2, FFR}
and about the degree and genus of curves in projective three-space \cite{MS2}. Here we use them
to prove two extensions of Reider's theorem on the positivity of adjoint linear series on a surface. We begin by recalling 
Reider's Theorem \cite{Re}.  

\medskip

Let $S$ be a smooth complex projective surface. 

\medskip

\nt {\bf Reider.}  Let $D$ be an ample (not necessarily effective) divisor on $S$. If 
\[ D^2 > 4 \ \mbox{and} \ D \cdot C > C^2 + 1 \]
for all curves $C \subset S$ satisfying $C^2 \le 0$, then $|K_S + D|$ is base point free, and if
\[ D^2 > 9 \ \mbox{and} \ D \cdot C > C^2 + 2 \]
for all $C \subset S$ with $C^2 \le 0$, then $|K_S + D|$ is very ample. 

\medskip

One extension due to Beltrametti, Francia and Sommese \cite{BFS}  is the following. If
\[ D^2 > (k+1)^2 \ \mbox{and} \ D \cdot C > C^2 + k \] 
for all $C \subset S$ with $C^2 \le 0$, then $\mbox{H}^1(S, \omega_S(D) \otimes I_Z) = 0$ for $Z \subset S$ of length $\le k$. 
This is a natural generalization of Reider's Theorem to the separation of $k$-jets. 

\medskip

Our first theorem is about the nef cone of the blow-up of $|K_S + D|^\vee$ along $S$.  

\medskip

\nt {\bf Theorem 7.1.} In the context of Reider's Theorem, suppose that 
\[ D^2 > 9 \ \mbox{and} \ D \cdot C > C^2 + 2 \ \mbox{when} \ C^2 \le 0 \] 
so that $S \subset |K_S + D|^\vee$ and let
$Y := \mbox{bl}_S(|K_S + D|^\vee)$ with exceptional divisor $E$.

\medskip

If in addition: 

\medskip

($a_1$) $D^2 \ge 16$ then $3H - E$ is nef on $Y$. 

\medskip

($a_2$) $D^2 > 16$ and $D \cdot C > C^2 + 3$ when $C^2 = -1,0$ then $3H - E$ is ample.  

\medskip

($b_1$) $D^2 \ge 25$ and $D\cdot C > C^2 + 3$ when $C^2 = -1,0$, then $2H - E$ is nef. 

\medskip

Note that $2H - E$ cannot be ample on $Y$ due to the presence of secant lines. However, if

\medskip

($b_2$) $D^2 > 25$ and $D \cdot C > C^2 + 4$ when $C^2= -1,0$, then $|2H - E|$ contracts: 
\[ (Y, \Sigma_1(S)) \rightarrow (Y_0, S^{[2]}) \] 
where $\Sigma_1(S)$ is the bundle of secant lines contracting onto the Hilbert scheme. 

\medskip

\nt {Remark.} The necessity of the additional  inequality in ($b_2$) is seen with the Veronese: 
\[ \PP^2 \subset \PP^5 \ \mbox{associated to} \ D = 5l  \] 
in which $2H - E$ on $Y$ is nef, but the secant variety is ``defective.'' Instead, in this case the contraction is a $\PP^2$-bundle over the dual projective plane
$(\PP^2)^\vee$. 

\medskip

For motivation, consider the case of a projective curve $C$. 

\medskip

\nt {\bf Riemann-Roch.}  If $\deg(D) > k$, then $|K_C + D|$ separates $k$-jets. 

\medskip

\nt This has extensions to the {\it equations} cutting out $C$, starting with: 

\medskip

\nt {\bf Castelnuovo.}  If $\deg(D) > 2$, the embedding of $C$ is projectively normal.

\medskip

Moreover, the homogeneous ideal $I_C$ is generated by quadrics and cubics. 

\medskip

\nt {\bf Saint-Donat.} If $\deg(D) > 3$, then $I_C$ is generated by quadrics. 

\medskip

Via the implications for a smooth subvariety $X \subset \PP^n$ and $Y = \mbox{bl}_X(\PP^n)$ 
\[ I_X \ \mbox{is generated by homogeneous polynomials of degree $\le d$} \] 
\[ \Downarrow \] 
\[ X \subset \PP^n \ \mbox{is scheme-theoretically cut out by polynomials of degree $d$} \] 
\[ \Downarrow \] 
\[ \ \mbox{the divisor} \ dH - E \ \mbox{on the blow-up} \ Y := \mbox{bl}_X(\PP^n) \ \mbox{is nef} \] 
\nt one can view Theorem 7.1 as a weaker version of the curve results, suggesting that the ideal
$I_S$ might be generated in low degree when the conditions on $D^2$ and $D \cdot C$ are met. 
What is missing is an understanding of the homogeneous rings: 
\[ \bigoplus \mbox{H}^0(Y, {\mathcal O}_Y(n(dH - E)))  \] 
and in particular whether or not they are generated in degree one. 

\medskip

As with Reider's Theorem, a curve $C$ for which $D \cdot C$ violates the inequality  in our theorem is directly 
responsible for the
failure of $dH - E$ to be nef (or ample). More interestingly, 
if there are no such curves but $dH - E$ still is not nef (or ample) then the failure is due to a special vector bundle 
on $S$ that comes into play with the failure of the inequality on $D^2$. On some surfaces, including  
$\PP^2, \PP^1 \times \PP^1$, $K3$ and abelian surfaces, there is detailed information about the (non)-existence of special vector bundles, 
which is consistent with sharper results.

\medskip

The proof follows a strategy outlined in Bayer-Macr\`i \cite{BM3}, in which the Bridgeland semi-stability of families of 
objects over a base (in this case $Y$) translates into the nefness of a {\it determinant} divisor class on $Y$. 
We begin by constructing a Serre family $C_{|K_S + D|^\vee}$ of objects
of the derived category $D^b(S)$ parametrized by $|K_S + D|^\vee$ which we pull back to the blow-up 
$Y$ where an elementary modification refines it to the {\it Drinfeld family} $B_Y$ that we use to implement the strategy. 
This is the higher dimensional analog of a construction of Drinfeld of a family 
of rank two vector bundles on a curve \cite{Be}.
The bulk of the work is then to find a chamber in the stability manifold of $S$ where the objects of $B_Y$ are 
Bridgeland semi-stable.

\medskip

Our second extension of Reider's Theorem concerns the Hilbert schemes $S^{[d]}$ 
of ideal sheaves of colength $d$ and may be seen as a generalization of the ampleness of $K_S + D$ when 
$D^2 > 9$ and $D \cdot C > C^2 + 2$ for all curves $C$ with $C^2 < 1$. 

\medskip

\nt {\bf Theorem 7.2.} For all surfaces $S$, the divisor classes 
\[ (K_S + D)^{[d]} - E\] 
are ample when
 $D$ is ample, $D^2 > 9d$ and $C \cdot D > C^2 + 2d$ for all curves with $C^2 < d$, 
where $(K_S + D)^{[d]}$ is the pull-back of the symmetrized divisor via the Hilbert-Chow morphism and
$E \subset S^{[d]}$ is the exceptional divisor of non-reduced schemes. 

\medskip

This, too, can be seen as the surface analogue of a well-known result for curves. 

\medskip

\nt {\bf Curves.} The Hilbert scheme is the symmetric product $C_d$ and divisor classes: 
\[ (K_C + D)_d - E \] 
are ample when $\deg(D) > 2d$. Indeed, in this case they are of the form: 
\[ (\deg(D) - 2d)x_d + 2\theta \] 
where $\theta$ is the pull-back via the Abel-Jacobi map of the (ample!) theta divisor on the 
Jacobian \cite{ACGH} and $x_d$ is a ``unit'' symmetrized divisor.  This shows
the bound is sharp. Ampleness when $\deg(D) = 2d$ is contingent on the Abel-Jacobi map 
being an embedding (or at least finite) which is in turn contingent on $C$ not being $d$-gonal. 

\medskip

Although there is no analogue of the Jacobian for $0$-cycles on a general surface by a theorem of Mumford \cite{Mu}, something similar 
seems to be happening relative to Theorem 7.2. If $f:S \rightarrow \PP^2$ is a finite map of degree $d$, then:
\[ D = f^*(-K_{\PP^2}) \ \mbox{satisfies} \ D^2 = 9d \]
and $(K_S + D)^{[d]} - E$ is not ample. As in the curve case, it is trivial when restricted to the embedded  $\PP^2 \subset S^{[d]}$.  
This is an interesting connection between the ampleness of such ``boundary'' divisors and  the {degree of irrationality} of the surface  (see \cite{BDELU}).

\medskip

There are many computations in the literature of the nef cones of Hilbert schemes for {special} surfaces  \cite{ABCH, Ry, BeC, BM1, BeC}, 
as well as some more recent, much more intricate computations of the larger movable cones  \cite{Ar, Bou, GrR}. 
The most general results on the nef cones of Hilbert schemes we are aware of are in \cite{BHLRS}, which also employs the Bayer-Macr\`i  strategy.
To our knowledge,  the results in this paper are the first  to apply to all surfaces and to condition ampleness on 
Reider-type inequalities. 

\medskip

\nt {\bf Acknowledgements.} The Drinfeld family of rank two vector bundles on a curve was explained to the first author by Robert Lazarsfeld. 
Extending it via additional blow-ups along secant varieties was the foundation of the author's PhD thesis \cite{Be}. Soon after, Michael Thaddeus utilized the
moduli spaces of Bradlow stable pairs as he replaced  blow-ups with Geometric Invariant Theory wall-crossing ``flips'' \cite{Th}. This innovation has
 flourished in the context of Bridgeland's stability conditions. The Bayer-Macr\`i strategy carried out in \S 7 
was partly inspired by joint work of the first author with Daniele Arcara, Izzet Coskun and Jack Huizenga \cite{ABCH}, in which it was noted that the walls for 
Bridgeland moduli of ideal sheaf objects on the projective plane were accompanied by continuously varying divisor classes on the Hilbert scheme. The first author would like 
to especially thank Daniele Arcara for sharing his obsession with Reider's Theorem (see e.g. \cite{AB2, LT}) and the conviction that it had to fit nicely with
Bridgeland stability conditions.

\newpage

\section{The Derived Drinfeld Construction}

We begin by reviewing the original Drinfeld construction (see \cite{Be}). 

\medskip

Let $|K_C + D|$ be an adjoint linear series on a curve $C$. Then via Serre duality, 
\[ |K_C + D|^\vee =  \PP(H^0(C,\omega_C(D)) = \PP(\Hom({\mathcal O}_C(-D), \omega_C)) = \Ext^1({\mathcal O}_C, {\mathcal O}_C(-D))^*/\CC^* \] 
is the projective space of lines through the origin in the vector space of {\bf extension} classes, which 
therefore generates non-split short exact sequences: 
\[ \epsilon: \ 0 \rightarrow {\mathcal O}_C(-D) \rightarrow V_\epsilon \rightarrow {\mathcal O}_C \rightarrow 0\ \mbox{for} \ \epsilon \in |K_C + D|^\vee  \] 

Via the {\it universal extension} over $C \times |K_C + D|^\vee$:
\[ 0 \rightarrow p^* {\mathcal O}_C(-D) \otimes q^*{\mathcal O}_{|K_C + D|^\vee}(1) \rightarrow V_{|K_C + D|^\vee} \rightarrow {\mathcal O}_{C\times |K_C + D|^\vee} 
\rightarrow 0 \] 
we obtain the {\it Serre} family 
$V_{|K_C + D|^\vee}$ of vector bundles
defining a rational map:
\[ \phi: |K_C + D|^\vee \dashrightarrow M_C(2, -D) \] 
to the moduli of {\bf semi-stable} rank two vector bundles $V$ on $C$ with $c_1(V) = -D$. 

\medskip

Now assume that
$C \subset |K_C + D|^\vee$ is an embedding, e.g. by taking $\deg(D) \ge 3$.
Then the points $p \in C$ parametrize the extensions that admit a (unique!) lift: 
\[ \begin{tikzcd}
    &&& {\mathcal O}_C(-p) \arrow[dl] \arrow[d] \\
    0 \arrow[r] & {\mathcal O}_C(-D) \arrow[r] & V_p \arrow[r] & {\mathcal O}_C \arrow[r] & 0
\end{tikzcd} \]
splitting the sequence after pulling back to the ideal sheaf $I_p = {\mathcal O}_C(-p)$. These are the most unstable of the bundles in the family, 
and the resulting exact sequence:
\[ 0 \rightarrow {\mathcal O}_C(-p) \rightarrow V_p \rightarrow {\mathcal O}_C(-D + p) \rightarrow 0 \] 
globalizes to a sequence centering the restriction $V_C := V_{|K_C + D|^\vee}|_{C \times C}$: 
\[ 0 \rightarrow {\mathcal O}_{C \times C}(-\Delta) \rightarrow V_C \rightarrow 
p^*{\mathcal O}_C(-D) \otimes q^* {\mathcal O}_C(K_C + D) \otimes {\mathcal O}_{C\times C}(\Delta) \rightarrow 0 \]  

The next step is to blow up $C \subset |K_C + D|^\vee$ to obtain: 
\[ Y := \bl_C(|K_C + D|^\vee) \ \mbox{with exceptional divisor} \ i:E \hookrightarrow Y \ \mbox{and} \ \pi:E \rightarrow C \] 
with fibers of $\pi$ isomorphic to $|K_C + D - 2p|^\vee$. Then the {\it pull-back} bundle:
\[ V_Y := \sigma^*V_{|K_C + D|^\vee}  \ \mbox{on} \ C \times Y \]  
restricts to a bundle $V_E$ with a quotient line bundle for the exceptional  {\bf divisor} $E$: 
\[ V_E = \pi^*V_C  \rightarrow L_E := \pi^* \left( {\mathcal O}_{C \times C}( (-D, K_C + D) + \Delta) \right) \] 
and kernel line bundle $\pi^* {\mathcal O}_{C \times C}(-\Delta)$. Then the elementary modification: 
\[ 0 \rightarrow W_Y \rightarrow V_Y \rightarrow i_*L_E \rightarrow 0 \]  
of $V_Y$ along $L_E$ is a new rank two vector bundle on $C \times Y$ such that: 

\medskip

(1) When $\epsilon \in Y - E$, then $W_\epsilon = V_\epsilon$ fits in 
\[ 0 \rightarrow {\mathcal O}_C(-D) \rightarrow W_\epsilon \rightarrow {\mathcal O}_C \rightarrow 0\] 
with no sub-bundle of degree $-1$ (or more). 

\medskip

(2) When $e \in E$ and $\pi(e) = p$, then $W_e$ fits in a non-split sequence: 
\[ 0 \rightarrow {\mathcal O}_C(-D + p) \rightarrow W_e \rightarrow {\mathcal O}_C(-p) \rightarrow 0 \] 
and so it, also, has no sub-bundle of degree $-1$ (or more). 

\medskip

The point is that via the elementary modification, the restriction $W_E$ satisfies: 
\[ 0 \rightarrow L_E(-E) \rightarrow W_E \rightarrow \pi^* {\mathcal O}_{C\times C}(-\Delta) \rightarrow 0 \] 
and then using ${\mathcal O}_E(-E) = {\mathcal O}_E(1)$, one inductively obtains (2).

\medskip

This generalizes to higher dimensions if we work in the bounded derived category $D^b(X)$ of coherent sheaves on 
a smooth projective variety $X$. 
Let $D$ be ample, and 
\[ X \subset |K_X + D| ^\vee  \] 
be an embedding of $X$ via the adjoint linear series (e.g. using Reider when $X = S$). 
Then the  Serre family $C_{|K_X + D|^\vee}$  is the {\bf cone} over the universal morphism: 
\[ \rightarrow C_{|K_X + D|^\vee}  \rightarrow {\mathcal O}_{X \times |K_X + D|^\vee} \stackrel u \rightarrow \left( p^*{\mathcal O}_X(-D) \otimes q^* {\mathcal O}_{|K_X + D|^\vee}(1)\right) [n]  \]  
in $D^b(X \times |K_X + D|^\vee)$ using Serre duality to identify: 
\[ |K_X + D|^\vee = \PP(H^0(X, \omega_X(D)) = \Hom_{D^b(X)} ({\mathcal O}_X, {\mathcal O}_X(-D)[n])^*/\CC^* \] 
as the space of lines through the origin in the vector space 
\[ \Ext^n_{{\mathcal O}_X}({\mathcal O}_X, {\mathcal O}_X(-D)) = \Hom_{D^b(X)} ({\mathcal O}_X, {\mathcal O}_X(-D)[n]) \] 

When $n = 1$, this recovers the Serre family of rank two vector bundles by rotating the distinguished 
triangle defining the cones $C_\epsilon$ once to the left, to obtain the previous extensions of line bundles. When $n > 1$, however, 
this rotation does not exhibit $C_{|K_X + D|^\vee}$ as a family of coherent sheaves, since it ``only'' yields a family of objects fitting in distinguished triangles:
\[  \rightarrow {\mathcal O}_X(-D)[n-1] \rightarrow C_\epsilon  \rightarrow {\mathcal O}_X \rightarrow  \] 

Nevertheless, we press on and identify the points $p\in X$ under the embedding. These are the cones for which 
there is an (again uniquely determined) lift: 
\[ \begin{tikzcd}
    && I_p \arrow[ld] \arrow[d] \\
    \arrow[r] & C_p \arrow[r] & {\mathcal O}_X \arrow[r] & {\mathcal O}_X(-D)[n] 
\end{tikzcd} \] 
and by an application of the octahedral axiom, we obtain distinguished triangles:  
\[ \rightarrow I_p \rightarrow C_p \rightarrow I_p^\vee (-D) [n-1] \rightarrow \] 
where 
$I_p^\vee = {\mathcal RHom}(I_p, {\mathcal O}_X)$
is the derived dual of $I_p$. This also globalizes  to: 
\[ \rightarrow {\mathcal I}_\Delta \rightarrow C_X \rightarrow {\mathcal I}^\vee_\Delta(-D,K_S + D) [n]  \rightarrow  \] 

We then let
$C_Y := L_\sigma^* C_{|K_X + D|^\vee}$ be the derived pull-back
and modify along the resulting morphism: 
\[ C_Y \rightarrow Ri_*\pi^* \left({\mathcal I}^\vee_\Delta(-D,K_S + D) [n]\right)  \] 
by taking the cone over the morphism to obtain a distinguished triangle: 
\[ B_Y \rightarrow C_Y \rightarrow Ri_*\pi^* \left({\mathcal I}^\vee_\Delta(-D,K_S + D) [n]\right) \] 
with the same breakdown of the elements of the family as in the curve case:

\medskip

(1) When $\epsilon \in Y - E$, then $B_\epsilon = C_\epsilon$ fits in 
${\mathcal O}_X(-D)[n-1] \rightarrow B_\epsilon \rightarrow {\mathcal O}_X$

\medskip

(2) When $e \in E$ and $\pi(e) = p$, then $B_e$ fits in 
${\mathcal I}_p^\vee(-D)[n-1] \rightarrow B_e \rightarrow {\mathcal I}_p $ 

\medskip

There is one important difference between curves and higher dimensions. 

\medskip

\nt {\bf Proposition 2.1.} If $\dim(X) > 1$, the objects $C_\epsilon$ are distinct and simple. 

\medskip

{\bf Proof.} The result follows directly from the distinguished triangles: 
\[ \begin{array}{cccccccccc} \cdots & {\mathcal O}_X(-D)[n-1] & \rightarrow & C_\epsilon & \rightarrow & {\mathcal O}_X \rightarrow \cdots \\
&&& \downarrow \\
\cdots & {\mathcal O}_X(-D)[n-1] & \rightarrow & C_{\epsilon'} & \rightarrow & {\mathcal O}_X \rightarrow \cdots 
\end{array} \] 
as well as  $\End({\mathcal O}_X) =  \CC \cdot \id = \End({\mathcal O}_X(-D)[n-1])$  and: 
\[ \Hom({\mathcal O}_X(-D)[n-1], {\mathcal O}_X) = 0  \] 

We can extend this to the Drinfeld family via: 

\medskip

\nt {\bf Proposition 2.2.} The objects $B_y$ of the family $B_Y$ are distinct and simple. 

\medskip

{\bf Proof.} Within a fiber of $\pi: E \rightarrow S$, the proof of Proposition 2.1 distinguishes 
$B_e$ from $B_{e'}$ since 
${\mathcal I}_p$ and ${\mathcal I}_p^\vee(-D)[n-1]$ are simple, and
\[ \Hom({\mathcal I}_p^\vee(-D)[n-1], {\mathcal I}_p) = 0 \] 
so that, additionally, each $B_e$ is simple. 

\medskip

This follows from the distinguished triangle ${\mathcal I}_p^\vee(-D)[n-1]$:  
\[ \rightarrow {\mathcal O}_S(-D)[n-1] \rightarrow {\mathcal I}_p^\vee(-D)[1] \rightarrow \CC_p \rightarrow \] 
and the fact that ${\mathcal I}_p$ is torsion-free. 

\medskip

Finally, $B_e$ knows the image of $e$ under the projection $\pi:E \rightarrow S$ from: 
\[ 0 \rightarrow \CC_p \rightarrow {\mathcal H}^0(B_e) \rightarrow {\mathcal I}_p \rightarrow 0 \] 
and does not coincide with any $C_\epsilon$ due to the torsion in ${\mathcal H}^0(B_e)$. \qed

\medskip

\nt Remark. This is a key aspect in which the Serre and Drinfeld families in higher dimensions 
differ from the families of rank two bundles. In the curve case, the bundles need not be distinct, or 
simple. Instead, it is the ``pair'' ${\mathcal O}_C(-D) \rightarrow E$  of the bundle together with a section of 
$E(D)$ that is simple, and stable for appropriate stability conditions. Since stable objects are always simple, it  is the simplicity of the objects
$C_\epsilon$ and $B_e$ that even make it plausible that they might be stable for some stability condition, as we 
will see next. 

\medskip

\section{Stability Conditions on a Surface} 

The Drinfeld family $W_Y$  ``refines''  the rational map defined by $V_{|K_C + D|^\vee}$: 
\[ \begin{tikzcd}
    Y \arrow[dr, dashed] \arrow[d] \\
    \ | K_C + D|^\vee \arrow[r, dashed, "\Phi"] & M_C(2, -D)
\end{tikzcd} \]
to the coarse moduli space of semi-stable rank two bundles $V$ with $c_1(V) = -D$. There is one stability of 
vector bundles on a curve $C$, based on the slope: 
\[ \mu(V) = \frac{\deg(V)}{\rk(V)} \] 
and a vector bundle $V$ is semi-stable if
$\mu(F) \le \mu(V)$ for all subbundles $F \subset V$.
When $X = S$, however, there are {\it families} of (Bridgeland) stability conditions. 

\medskip

We choose an ample polarization $H$ on $S$, which we normalize so that $H^2 = 1$.  When 
$D$ is an ample {\it integral} divisor, we achieve this by setting
$H = D/||D||$ and 
then there are two $H$-slope analogues of the slope of vector bundle on a curve. 

\medskip

\nt {\bf Mumford Slope} is defined on the {\it torsion-free} sheaves $F$:
\[ \mu_H(F) = \frac{c_1(F) \cdot H}{\rk(F)} \] 

\nt {\bf Gieseker/Simpson Slope} is defined on {\it pure dimension one} torsion sheaves $\tau$:
\[ \nu_H(\tau) = \frac{\ch_2(\tau)}{c_1(\tau) \cdot H} \] 

To hold on to more information, one uses instead the central charges: 
\[ Z_{\mu}(F) = (-c_1(F) \cdot H, \rk(F)), \ \ Z_\nu(\tau) = ( -\ch_2(\tau), c_1(\tau) \cdot H)    \] 
underlying the Mumford and Gieseker slopes. 
These are linear on the $K$-group  and non-zero for non-zero sheaves $F$ (and $\tau$), with
$0 < \arg(Z({\mathcal E})) < 1$,  where $\arg(re^{\pi i \theta}) := \theta$  
for sheaves ${\mathcal E} = F$ or $\tau$.  

\medskip

Since arg$(Z)$ imposes the same ordering as the slope, one makes the following
definition, enlarging the categories to include sheaves with $\arg(Z({\mathcal E})) = 1$. 

\medskip

\nt {\bf Definition.} An arbitrary  torsion sheaf $\tau$ on $S$ is $H$-semistable if: 
\[  \arg(Z_\nu(\tau)) \le \arg(Z_\nu(\tau')) \ \mbox{for all quotients} \ \tau \rightarrow \tau' \] 
A sheaf ${\mathcal F}$ with no torsion in dimension zero is $H$-semistable if: 
\[  \arg(Z_\mu(\mathcal F)) \le \arg(Z_\mu(\mathcal F')) \ \mbox{for all quotients} \ \mathcal F \rightarrow \mathcal F'  \] 
and $H$-stable if the inequality is always strict and strictly $H$-semistable otherwise. 
 
 \medskip
 
 \nt Note. A subsheaf of a coherent sheaf with no torsion in dimension zero also has no torsion in dimension zero, but to get a quotient
 with the same property, one needs to saturate the subsheaf. This is the reason we switched to quotients.

\medskip

The good properties of these slopes include: 

\medskip

\nt {\bf Schur's Lemma.} Stable sheaves  are simple, i.e. End$(E) =  \CC \cdot \id_E$, and
\[ \Hom(E_1, E_2) = 0 \ \mbox{if $E_1,E_2$ are semi-stable and} \ \arg(Z(E_1)) > \arg(Z(E_2)) \]

\nt {\bf Finiteness of Filtrations.} 

\medskip

(a) Semi-stables have finite Jordan-Holder filtrations with stable sub-quotients. 

\medskip

\nt Note. In order to define this as an increasing filtration, one saturates the inclusions. 

\medskip

(b) Arbitrary sheaves (of each given type) have unique finite Harder-Narasimhan filtrations with semi-stable 
sub-quotients with strictly decreasing slopes. 

\medskip

\nt {\bf Projective Moduli.} Stability is an open condition on flat families, semi-stability is also open with a replacement property and 
projective coarse moduli 
spaces
\[ {\mathcal M}_S(r, c_1, d) \ \mbox{and} \ {\mathcal M}_S(0,D,d) \] 
exist for each set of Chern character invariants of torsion-free/torsion sheaves. 

\medskip

\nt {\bf Bogomolov Inequality.} Mumford semi-stable torsion-free sheaves satisfy
\[ c_1(F)^2 \ge 2 \cdot \ch_2(F) \cdot \rk(F) \] 
from which we get the {\bf $H$-Bogomolov inequality}, via  the Hodge index thoerem:
\[ (c_1(F) \cdot H)^2 \ge 2 \cdot \ch_2(F) \cdot \rk(F) \] 
with equality if and only if the first is an equality and $c_1(F) = \lambda H$ for some $\lambda$. 

\medskip

There are various ways to introduce geometric Bridgeland stability conditions on a polarized surface 
as deformations of  either of the ``classical'' stability conditions. We choose here to deform the Gieseker slope, i.e.
the associated central charge:
\[ Z_\nu(\tau) = (-\ch_2(\tau), \ c_1(\tau) \cdot H) \] 
by introducing the rank of an {\bf arbitrary} object $E$ of $D^b(S)$ via: 
\[ Z_{(x,y)}(E) :=  ( -\ch_2(E) + x\cdot \rk(E), \ c_1(E) \cdot H + y\cdot \rk(E))  \] 
or, equivalently, as a family of slopes:
\[ \nu_{(x,y)}(E) = \frac{\ch_2(E) - x\cdot \rk(E)}{c_1(E)\cdot H + y\cdot \rk(E)} \] 
indexed by points of the $(x,y)$-plane. 

\medskip

A preliminary observation of these slopes gives:

\medskip

(i) the slopes of sheaves $\tau$ of pure dimension one do not vary with $(x,y)$.   

\medskip

(ii) $\arg(Z_{(x,y)}(\CC_p)) = 1$ is maximal for all skyscraper sheaves $\CC_p$ and all $(x,y)$. 

\medskip

(iii) Each object $E$ of $D^b(S)$ with $\rk(E) \ne 0$ has a unique  {\bf indeterminate point} 
\[ I(F) = \left( \frac{\ch_2(F)}{\rk(F)}, - \frac{c_1(F)\cdot H}{\rk(F)} \right) = \left( \frac{\ch_2(F)}{\rk(F)}, - \mu(F) \right) \] 
defined by the property that $Z_{I(F)}(F) = (0,0)$. Moreover, 

\medskip

\nt {\bf Corollary} ($H$-Bogomolov) The point of  indeterminacy $(x,y) = I(F)$ of a 
Mumford semi-stable torsion-free sheaf $F$ always satisfies $x \le \frac {y^2} 2$.

\medskip

We want to say that $Z_{(x,y)}$ defines a stability condition when, conversely, 
$x > \frac {y^2} 2$ 
and that in that case, we can constrain $\arg(Z_{(x,y)}(E)) \in (0,1]$ just 
as the arg's of coherent sheaves on a curve were constrained. 
This is the case for all torsion sheaves by the observation above, but not for torsion-free sheaves. However, 
Bridgeland first explained how a family of tilts of the category of coherent sheaves to 
``coherent-adjacent'' abelian categories achieves this. 

\medskip

For an {\bf arbitrary} coherent sheaf ${\mathcal F}$ on $S$, let $\tau({\mathcal F}) \subset {\mathcal F}$ be the torsion subsheaf and  
consider the extended Harder-Narasimhan filtration by saturated subsheaves:
\[ \tau({\mathcal F})  =: {\mathcal F}_0 \subset {\mathcal F}_1 \subset \cdots \subset {\mathcal F}_l \subset {\mathcal F} \] 
chosen so that
${\mathcal F}_{i}/{\mathcal F}_{i-1}$ are Mumford semi-stable, of strictly decreasing slopes. 

\medskip

\nt {\bf Definition 3.2.} For any non-zero coherent sheaf ${\mathcal F}$ on $S$, let:
\[ \mu_{\min}({\mathcal F}) = \left\{ \begin{array}{l} +\infty \ \mbox{if} \ {\mathcal F} = \tau({\mathcal F}) \\ 
\mu({\mathcal F}/{\mathcal F}_l) \ \mbox{otherwise} \\ 
\end{array} \right. \ \mbox{and} \ \ 
\mu_{\max}({\mathcal F}) = \left\{ \begin{array}{l} +\infty \ \mbox{if} \ \tau({\mathcal F}) \ne 0 \ \mbox{and}  \\  
\mu({\mathcal F}_1)  \ \mbox{otherwise} \\ 
\end{array} \right. \]

The following is a direct consequence of the Harder-Narasimhan filtrations: 

\medskip

\nt {\bf Lemma 3.3.} For $y \in \RR$, the pair $({\mathfrak F}_y, {\mathfrak T}_y)$ of full subcategories of Coh$(S)$ with: 
\[ ob({\mathfrak T}_y) =  \{ \mbox{coherent sheaves ${\mathcal F}$ on $S$ with}  \ \mu_{\min}({\mathcal F}) > y \} \] 
\[ ob({\mathfrak F}_y) =  \{ \mbox{coherent sheaves ${\mathcal F}$ on $S$ with}  \ \mu_{\max}({\mathcal F}) \le y \} \] 
is a {\bf torsion pair}, i.e. $\Hom(T,F) = 0$ for $T \in {\mathfrak T}_y$ and $F \in {\mathfrak F}_y$ and each ${\mathcal F}$ sits in: 
\[ 0 \rightarrow T_y({\mathcal F}) \rightarrow {\mathcal F} \rightarrow F_y({\mathcal F}) \rightarrow 0 \] 
a uniquely determined short exact sequence. 

\medskip

\nt {\bf Definition 3.4.} For each $y \in \RR$, the {\bf Mumford tilt} at $y$ is the category: 
\[ {\mathcal A}_y \ \mbox{equipped with torsion pair}  \ ({\mathfrak T}_y, {\mathfrak F}_y[1]) \]
consisting of all objects $E$  with cohomology in degrees $-1$ and $0$ only, and: 
\[ {\mathcal H}^{-1}(E) \in \mbox{ob}({\mathfrak F}_y) \ \mbox{and} \ {\mathcal H}^0(E) \in \mbox{ob}({\mathfrak T}_y) \] 

These are abelian categories  (hearts of an underlying $t$-structure) rigged so that:
\[ \arg(Z_{(x,y)}(E)) \in (0,1] \ \mbox{for} \  E \in \mbox{ob}({\mathcal A}_{-y}) \] 
and it is a deep theorem 
due to Bridgeland in the K3 case and Bayer-Macr\`i and Toda in general that these 
are good slopes from the point of view of stability. 

\medskip

\nt {\bf Theorem} (see \cite{MS1}) For $(x,y)$ in the {\it stability region} $x > \frac 12y^2$, the pairs:
\[ (Z_{(x,y)}, {\mathcal A}_{-y}) \] 
are a continuous family of Bridgeland stability conditions. In particular, 

\medskip

(a) $Z_{(x,y)}$-semi-stable objects of ${\mathcal A}_{-y}$ all have finite Jordan-Holder filtrations

\medskip

(b) Arbitrary objects of ${\mathcal A}_{-y}$ have finite $Z_{(x,y)}$-Harder-Narasimhan filtrations

\medskip

(c) All $(x,y)$-semi-stable objects $E \in {\mathcal A}_{-y}$ satisfy the $H$-Bogomolov inequality. 

\medskip

(d) For fixed invariants $(r, D, d)$, there is a locally finite chamber decomposition 
of the stability region $x > \frac 12y^2$ with linear walls such
that the moduli of semi-stable objects do not change as $(x,y)$ varies within a chamber.  

\medskip

The walls and chambers are features of the {\it variation} of the stability conditions.
For now, we focus on some preliminary examples. 

\medskip

\nt {\bf Example 1.} Skyscraper sheaves are stable. 
Since $\CC_p \in {\mathcal A}_{-y}$ has maximal phase:
\[ \arg(Z_{(x,y)}(\CC_p)) = 1 \] 
 it follows that skyscraper sheaves are semi-stable for all $(x,y)$. When $y \not\in \QQ$, the finite-length sheaves are the only 
objects of maximal phase. When $y \in \QQ$, there are also the shifts
$F[1]$
of a Gieseker semi-stable torsion-free sheaf of slope $-y$. But
$\Hom(F[1], \CC_p) = 0$ and it follows that the $\CC_p$ are stable. 

\medskip

\nt {\bf Example 2.} Let $F$ be a Mumford-semi-stable torsion-free sheaf of slope $\mu(F)$. 

\medskip

$\bullet$ $F[1]$ is $(x,-\mu(F))$-semistable for all $x > \frac 12\mu(F)^2$. 

\medskip

$\bullet$ $F[1]$ is $(x,-\mu(F))$-stable when $F = V$ is locally free and Mumford stable. 

\medskip

As in Example 1, the stability of $V[1]$ in the second bullet point follows from:
\[ \Hom(\CC_p, V[1]) = \Ext^1(\CC_p,V) = \Ext^1(V, \CC_p \otimes \omega_S)^\vee = H^1(S, \CC_p \otimes \omega_S \otimes V^*)^\vee = 0 \] 

Otherwise, $F[1]$ is either strictly $(x,-\mu(F))$-semistable by virtue of having a Mumford-stable subsheaf of the same slope, {\bf or} by 
virtue of the rotation:
\[ 0 \rightarrow T \rightarrow F[1] \rightarrow F^{**}[1] \rightarrow 0  \ \mbox{of the sequence} \  0 \rightarrow F \rightarrow F^{**} \rightarrow T \rightarrow 0 \]  
Conversely, an extension:
\[ \epsilon: \ 0 \rightarrow V[1] \rightarrow E \rightarrow T \rightarrow 0 \ \ \mbox{for} \ \epsilon \in \Ext^ 2(T,V) = \Hom(V,T \otimes \omega_S)^\vee \] 
also produces an $(x,y)$-semistable object $E$. 
If it is  stable, it is the shifted derived dual  $F^\vee[1]$ of a stable torsion-free sheaf $F$ with ordinary dual $F^* = V$. 

\medskip

It follows from openness of stability that a Mumford stable vector bundle $V$ is 
$(x,y)$-stable in a neighborhood of each $(x, -\mu(V))$ 
but $V$ needs to be interpreted as a stable object of $D^b(S)$, i.e. as the  {\it shift} of a stable object in ${\mathcal A}_{-y}$, since $V \in {\mathcal A}_{-y}$ for $y > -\mu(V)$, but 
$V[1] \in {\mathcal A}_{-y}$ for $y \le -\mu(V)$. Also
\[ 0 \rightarrow F \rightarrow V = F^{**} \rightarrow T \rightarrow 0 \] 
shows that a Mumford-stable $F$ is $(x,y)$-stable in a half-neighborhood 
$y > -\mu(F)$ of each point $(x,-\mu(F))$. The following  is a consequence of the Theorem of Bayer-Macr\'i-Toda
and the nature of the {\it destabilizing walls} for objects $F \in {\mathcal M}_{(x,y)}(r,c_1,d)$, pointing out 
the criticality of  (c) and (d) of the Theorem to understanding the walls and chambers for moduli.

\medskip

$\bullet$ If  $r > 0$ and $c_1^2 \ge 2rd$, the coarse moduli space of Bridgeland-semistable objects ${\mathcal M}_{(x,y)}(r,c_1,d)$ 
coincides with the (projective) coarse  moduli ${\mathcal M}_S^G(r,c_1,d)$ of  {\it Gieseker} semistable torsion-free sheaves for all $(x,y)$ in a {\it sector}: 
\[\left(  \left( \frac dr, \frac {c_1\cdot H}r  \right) + z \ | \ 0 < \arg(z) < \theta^+ \right) \cap \left( x > \frac 12y^2 \right) \]

(b) The moduli ${\mathcal M}_{(x,y)}(-r, -c_1, -d)$ coincide with 
${\mathcal M}^{G}_S(r,-c_1,d)$  in a sector: 
\[ \left(  \left( \frac dr, \frac{c_1\cdot H}r  \right) + z \ | \ \theta^- < \arg(z)  < 0 \right) \cap \left( x > \frac 12y^2 \right) \]
via the map $F \mapsto F^\vee[1]$. 

\medskip

There are several interesting features here, including the appearance of the new moduli spaces of 
Gieseker semi-stable torsion-free sheaves. This is a refinement of the moduli of Mumford semi-stable 
sheaves turning strictly semi-stable sheaves into Gieseker stable, strictly semistable or unstable 
sheaves depending upon the 
second chern character (see \cite{BeM}). This is quite interesting, but we will not need it here.
Another  interesting and very relevant feature is the sector shape of the chamber. This is 
due to the nature of the {\bf walls} for objects with nonzero rank.

\medskip

\nt {\bf Proposition 3.5.} Let $E \subset F$ be a subobject in ${\mathcal A}_{-y}$
of a Mumford-stable torsion-free sheaf $F$ and suppose
$E$ and $F$ are $(x,y)$-semistable with the same $(x,y)$-slope. Then:

\medskip

(a) $E= {\mathcal H}^0(E)$ is a  sheaf, i.e. an element of ${\mathfrak T}_{-y}$ by virtue of the sequence: 
\[ 0 \rightarrow {\mathcal H}^{-1}(E) \rightarrow 0 \rightarrow {\mathcal H}^{-1}(F/E) \rightarrow {\mathcal H}^0(E) \rightarrow F \rightarrow {\mathcal H}^0(F/E) \rightarrow 0 \] 
on the cohomology sheaves of the objects $E,F,F/E$ in the category ${\mathcal A}_{-y}$. 

\medskip

(b) The line through $I(E)$ and $I(F)$ passes through $(x,y)$. Moreover, $I(E)$ lies 
above $I(F)$ and between it and the parabola (see Figure 1). 

\medskip

\begin{figure}
\begin{center}
    \begin{tikzpicture}
        \begin{axis}[
            axis lines = center,
            xmin = -5, xmax = 15, xticklabels = {,,},
            ymin = -5, ymax = 5, yticklabels = {,,},
            unit vector ratio = 1 1
        ]
        \addplot[color=red, samples=100, domain=0:15]{sqrt(2*x)};
        \addplot[color=red, samples=100, domain=0:15]{-sqrt(2*x)};

        \addplot[color=OliveGreen, domain=-5:15]{x/3+1/6};
        \draw[blue, fill] (axis cs: -2, -0.5) circle (2pt) node [below] {$I(F)$};
        \draw[Purple, fill] (axis cs: -0.5, 0) circle (2pt) node [above] {$I(E)$};
        \draw[black, fill] (axis cs: 5.5, 2) circle (2pt) node [below right] {$(x, y)$};
        \end{axis}
    \end{tikzpicture}
\end{center}
\caption{When $E \subset F$ causes $F$ to fail to be  $(x,y)$ stable.}

\end{figure}
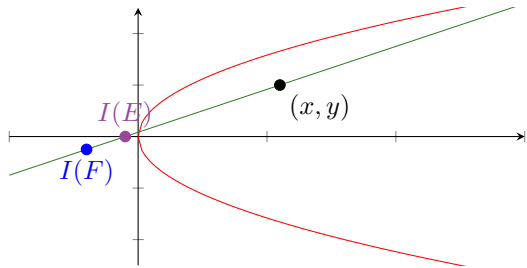

{\bf Proof.} By (c) of the Theorem, each of $I(E)$ and $I(F)$ lie outside the parabola. Moreover, $Z_{(x,y)}(E)$ and $Z_{(x,y)}(F)$ are on the same (real) line in $\CC$, from 
which it follows that $(x,y)$ (in the ``phase'' space) is on the line through $I(E)$ and $I(F)$. Finally, we've seen $(x,y) = I(F) + z$ with $\arg(z) > 0$. Coupling that with 
$E \in {\mathfrak T}_{-y}$ for some
$y > -\mu(E)$, it follows that $I(E)$ is between $I(F)$ and the parabola. 

\medskip

\nt Note. If $E/F$ has non-zero rank, $I(E/F)$ is on the opposite side of the parabola. 

\medskip

\nt {\bf Applications.} (a) ${\mathcal O}_S$  is $(x,y)$ stable for all $y > 0$ and $(x,y)$ in the parabola. 

\medskip

More generally, if $F$ is $H$-stable and the $H$-Bogomolov inequality is an equality, 
then $F$ is $(x,y)$ stable for all 
$y > -\mu(F)$ and $x > \frac 12y^2$. Similarly, in that case $F^\vee[1]$ is $(x,y)$ stable for 
all $y < \mu(F)$. The reason for this is simple. There is no room to place $I(E)$ between $I(F)$ and the parabola!

\medskip

(b) Consider the ideal sheaf ${\mathcal I}_p$ of a point $p \in S$ with $I({\mathcal I}_p) = (-1,0)$. 
Since the line through $I({\mathcal I}_p)$ tangent to the parabola meets the parabola at $(1,\sqrt 2)$, 
there is room for destabilizing subobjects of ${\mathcal I}_p$. 
For example, consider ${\mathcal O}_S(-C) \subset {\mathcal I}_p$ corresponding to a curve $C$ through $p$. Then (see Figure 2)
\[ I({\mathcal O}_S(-C)) = \left( \frac{C^2}2, C \cdot H \right) \] 
lies between $(-1,0)$ and the parabola (and therefore destabilizes $I_p$) when: 

\medskip

$\bullet$ $C^2 = -1$ and $C \cdot H < \frac 1{\sqrt 8}$ or $D^2 > 8(C\cdot D)^2$ for an integral ample $D = H \cdot ||D||$.

\medskip

$\bullet$ $C^2 = 0$ and $C \cdot H < \frac 1{\sqrt 2}$ or $D^2 > 2(C \cdot D)^2$ for $D$. 

\medskip

$\bullet$ $C^2 = 1$ and $1 \le C\cdot H < \frac 3{\sqrt 8}$ or $(C \cdot D)^2  \ge D^2  > \frac 89 (C\cdot D)^2$ for $D$. 

\medskip

\begin{figure}
\begin{center}
    \begin{tikzpicture}
        \begin{axis}[
            axis lines = center,
            xmin = -3, xmax = 7, xticklabels = {,,},
            ymin = -2, ymax = 3, yticklabels = {,,},
            unit vector ratio = 1 1
        ]
        \addplot[name path = para1, color=red, samples=100, domain=0:7]{sqrt(2*x)};
        \addplot[name path = para2, color=red, samples=100, domain=0:7]{-sqrt(2*x)};
        \addplot[name path = floor, draw=none, domain=-1:0]{0};

        \addplot[name path = line, color=OliveGreen, domain=-3:7]{x/sqrt(2)+1/sqrt(2)};
        \addplot[color=gray] fill between[of = line and floor, soft clip = {domain= -1:0}];
        \addplot[color=gray] fill between[of = line and para1, soft clip = {domain = 0:1}];
        \draw[black, fill] (axis cs: -1, 0) circle (2pt) node [above left] {$I(I_p)$};
        \end{axis}
    \end{tikzpicture}
\end{center}

\caption{The region for $I({\mathcal O}_S(-C))$ when $p\in C$ is destabilizing}
\end{figure}
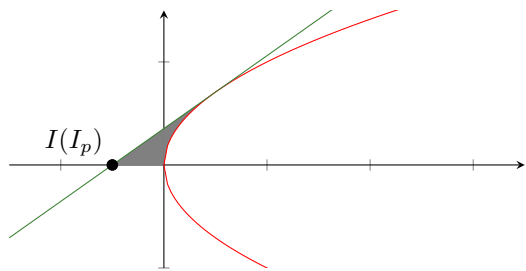

\nt {\bf Important Remark.} The stability region of $I_p$ depends on $H$. When $\Pic(S) = \ZZ$, there is only one 
choice for $H$ (and there are no curves with $C^2 = -1$ or $0$) but when the Picard rank is larger than one, the choice matters. 
From this point of view, let's think more about the inequalities $D^2 > n (C \cdot D)^2$. 

\medskip

$\bullet$ $C^2 = -1$. In this case, we can achieve $C \cdot D = 1$ with large $D^2$ by letting: 
\[ \sigma: S \rightarrow S_0 \ \mbox{be the blow-down of an exceptional curve $C = \PP^1$ and}  \ D = \sigma^*D_0 - E \] 
for $D_0$ sufficiently ample on $S_0$. Thus, the condition $D^2 > 8(D \cdot C)^2$ should be thought of as saying that 
$H$ is near enough to the pull-back of an ample $H_0$. 
Interestingly, the blow-down is itself a wall-crossing of moduli spaces, since:  
\[ 0 \rightarrow {\mathcal O}_S(-C) \rightarrow I_p \rightarrow {\mathcal O}_C(-p) \rightarrow 0 \] 
are short exact sequences parametrized by: 
\[ \Ext^1({\mathcal O}_C(-p)), {\mathcal O}_S(-C)) = \Ext^1({\mathcal O}_S(-C), {\mathcal O}_C(-p) \otimes \omega_S)^\vee = \mbox{H}^1(C,\omega_C(-p)) ^\vee \] 
by the adjunction formula. But $C = \PP^1$, so $\mbox{H}^1(C, \omega_C(-p)) = \mbox{H}^1(\PP^1, {\mathcal O}_{\PP^1}(-3))$ projectivizes to $\PP^1$, corresponding to the ideal sheaves $I_p$ for $p \in C = \PP^1$.

\medskip

On the other hand, the extensions in the opposite direction (after wall-crossing): 
\[ 0 \rightarrow  {\mathcal O}_C(-p) \rightarrow E \rightarrow  {\mathcal O}_S(-C) \rightarrow 0 \] 
are indexed by
$ \Ext^1({\mathcal O}_S(-C), {\mathcal O}_C(-p)) = \mbox{H}^1(C, {\mathcal O}_C(-p) \otimes {\mathcal O}_S(C)) = \mbox{H}^1(\PP^1, {\mathcal O}_{\PP^1}(-2))$ 
which projectivizes to a {\bf point}. From this we conclude that
the moduli space $S$ of ideal sheaves $I_p$ contracts to $S_0$ at the wall (and beyond). 

\medskip

$\bullet$ $C^2 = 0$. We get a similar picture when $C$ is the fiber of a ruled surface $S \rightarrow B$, 
except for the fact that the opposite direction extensions: 
\[ 0 \rightarrow  {\mathcal O}_C(-p) \rightarrow E \rightarrow  {\mathcal O}_S(-C) \rightarrow 0 \] 
all split, since they are parametrized by $\mbox{H}^1(\PP^1, {\mathcal O}_{\PP^1}(-1)) = 0$. This means 
that at the wall, the moduli space contracts to $B$, but does not extend further.  

\medskip

$\bullet$ $C^2 = 1$. For  $\PP^2$, it is ${\mathcal O}_{\PP^2}(-1)^2 \rightarrow {\mathcal I}_p$ that gives a non-trivial wall.  

\medskip

\nt {\bf Second Remark.} This does not exhaust all the walls that one might encounter. 
For example, there may be ``higher rank'' walls corresponding to: 
\[ E \rightarrow I_p \] 
where $E$ (for example) is a Mumford semi-stable bundle with a nonzero map to $I_p$ corresponding to a section of the dual bundle
vanishing at $p$. This map is not an injective map of coherent sheaves, but it may be injective 
in a category ${\mathcal A}_{-y}$ for a stability condition $(x,y)$ for which $I(E)$ lies between $(-1,0)$ and the parabola.

\medskip

We will deal with such hypothetical walls  in \S 5. 

\medskip

\section{Stability Conditions and the Serre Family} 

We start by finding a chamber in the stability manifold where the Serre family
$C_{|K_S + D|^\vee}$
is a family of $(x,y)$-{\bf stable} objects. Their class in the $K$-group is: 
\[ [{\mathcal O}_S] - [{\mathcal O}_S(-D)] \ \mbox{with Chern character} \ \left(0, D, -\frac{D^2}2 \right) \] 

There is no point of indeterminacy for such objects. Instead, if: 
\[ 0 \rightarrow F \rightarrow C_\epsilon \rightarrow E \rightarrow 0 \] 
is a short exact sequence of $(x,y)$-semistable objects in ${\mathcal A}_{-y}$, then: 
\[ I(E) - I(F) = \left( \frac{D^2}{2\rk(F)}, \frac{D \cdot H}{\rk(F)}\right) = \frac 1{\rk(F)} \left( \frac{D^2}{2}, {||D||} \right)  \] 
i.e. all the walls are {\bf parallel}, with the same slope
$2/{||D||}$.

\medskip 

\nt {\bf The Genesis Wall $W_0$.} This wall, extending between points {\bf on} the parabola: 
\[ (0,0) = I({\mathcal O}_S) \ \mbox{to} \ \left( \frac{D^2}2, ||D||  \right) = I({\mathcal O}_S(-D)[1])  \] 
{\it generates} the family $C_{|K_S + D|^\vee}$ in the sense that: 

\medskip

(a) Each $C_\epsilon$ is strictly semi-stable for $(x,y) \in W_0$, 

\medskip

(b) Each $C_\epsilon$ is {\bf stable} for $(x,y)$ between $W_0$ and a next wall $W_1$.  

\medskip

This follows again from (c) and (d) of the Theorem of Bayer-Mac\`i-Toda and the fact that 
${\mathcal O}_S$ and ${\mathcal O}_S(-D)[1]$ are both $(x,y)$-stable in the entire region
$0 < y \le ||D||$. 

\medskip

\nt {\bf The Next Wall $W_1$} crosses the $x$-axis at $(a,0)$ for $a \in [-1,0)$ and contains either: 

\medskip

(i) $I(F)$ for some torsion-free $F \subset C_\epsilon$ of rank $\ge 2$ and some $\epsilon \in |K_S + D|^\vee$. 

\medskip

(ii) $I({\mathcal I}_Z(-C))$ for a curve $C \subset S$ with ${\mathcal I}_Z(-C) \subset C_\epsilon$ for some $\epsilon$. 
 
 \medskip
 
 (iii) $I({\mathcal I}_p) = (-1,0)$ which we know satisfies ${\mathcal I}_p \subset C_p$ for $p\in S \subset |K_S + D|^\vee$. 

\medskip

By definition, the wall $W_1$ passes through $I(E)$ for some subobject $E \subset C_\epsilon$ with:
\[ \arg(Z_{(x,y)}(E)) = \arg(Z_{(x,y)}(C_\epsilon)) \] 

The cohomology sequence  then gives
$0 \rightarrow {\mathcal H}^{-1}(E) \rightarrow {\mathcal O}_C(-D) \rightarrow {\mathcal H}^{-1}(C_\epsilon/E)$
and since each of these is torsion-free, it follows that either: 
\[ {\mathcal H}^{-1}(E) = 0 \ \mbox{(and $E$ is a coherent sheaf) or else} \ {\mathcal H}^{-1}(E) = {\mathcal O}_S(-D) \] 

But in the latter case, we have a diagram: 
\[ \begin{array}{ccccccccccc} 0 & \rightarrow &  {\mathcal O}_S(-D)[1] & \rightarrow & E & \rightarrow & {\mathcal H}^0(E) & \rightarrow & 0 \\
&& || && \downarrow && \downarrow \\
0 & \rightarrow &  {\mathcal O}_S(-D)[1]  & \rightarrow & C_\epsilon & \rightarrow & {\mathcal O}_S & \rightarrow & 0 \end{array} \] 
from which we conclude $\arg(Z_{(x,y)}({\mathcal H}^0(E))) \ge \arg(Z_{(x,y)}({\mathcal O}_S)))$, contradicting the stability of ${\mathcal O}_S$.  
We now get to the crux of the matter. 

\medskip

If
 $F \subset  C_\epsilon$  
witnesses the first wall with $E = C_\epsilon/F$, then $F$ is torsion-free due to: 
 \[ \rightarrow {\mathcal H}^{-1}(E) \rightarrow F \rightarrow {\mathcal O}_S \rightarrow \] 
 from the cohomology sequence. 
 
 \medskip
 
 Moreover, $I(F)$ and $I(E)$ are  on {\bf opposite sides} of the parabola, with: 
\[ \mu_H(F) = \frac{c_1(F)\cdot H}{\rk(F)} \ \mbox{and} \ \mu_H(E) = \frac{(D - c_1(F)) \cdot H}{\rk(F)}  \] 

In particular, their difference  is 
${ ||D||}/{\rk(F)}$ 
which is therefore an {\it upper bound} for  the {\bf height} of the wall $W_1$. On the other hand, we can easily calculate this height in terms of $a$, since it 
is obtained by intersecting: 
\[ \left(y = \frac 2{||D||}x - a \frac 2{||D||}\right) \ \mbox{with} \ \left(x = \frac 12 y^2 \right) \]   

This becomes $y^2 - ||D|| y - 2a = 0$
and the height is the difference of the roots. Thus,  
\[ \sqrt{D^2 + 8a} \le \frac{||D||}r \] 
simplifying to
\[ D^2 \le (-8a) \left( \frac {r^2}{r^2 - 1}\right)  \] 
is required for the first wall $W_1$ to be witnessed by an $F$ of rank $r > 1$. 

\medskip

In particular, if:
\[ D^2 > \frac{32}{3}  \ \mbox{or, more succinctly,} \ D^2 > 10 \] 
then $W_1$ is witnessed by a rank one sheaf ${\mathcal I}_Z(-C) \subset C_\epsilon$ for a curve $C \subset S$ and subscheme $Z \subset S$ of length $l$. The
line through $I({\mathcal I}_Z(-C)) = \left( \frac {C^2}2 - l, - C \cdot H \right)$ with slope $ \frac 2{||D||}$
meets the $x$-axis at $(a,0)$ for  
\[ -2a - 2l = C \cdot D - C^2 \]  
so if $C \cdot D > C^2 + 2$ for all $C$, then $a < -1$ and $W_1$ is witnessed by $I_p \subset C_p$. 

\medskip

This very nearly matches Reider's condition for the linear series to be very ample. One needs to know that it suffices to check curves with $C^2 \le 0$. From:
\[ \mu({\mathcal O}_S(-C)) = -C \cdot H   \ \mbox{and} \ \mu(C_\epsilon/{\mathcal O}_S(-C)) = (C - D)\cdot H \] 
and the fact that the points of indeterminacy are on opposite sides of the parabola of stability conditions,  it suffices to check $C$ for which: 
\[ C\cdot D < \frac 12 D \cdot D \] 

Together with the Hodge Index Theorem this implies 
$C^2 < \frac 12 C \cdot D$
and then the inequality $C\cdot D > C^2 + 2$ follows for $C^2 > 0$ unless $C^2 = 1$ and $C \cdot D = 3$. But when we also make 
use of $D^2 > 10$ (or just $9$), we get 
$9 = (C \cdot D)^2 \ge D^2 \cdot C^2  > 10$ and a contradiction to that case.   

\medskip

\nt {\bf Remark.} If we could rule out rank two witnesses in some other way, we'd shave down to $D^2 > 9$ and perfectly
match the Reider criterion for very ampleness. This seems more likely to be the correct inequality. See Martinez (\cite{Ma}).

\section{Some Estimates Regarding Stability Regions for the Sheaves ${\mathcal I}_Z$} 

Let $S$ be a surface with an ample $H$  normalized so that $H^2 = 1$. 

\medskip

Let $S^{[n]}$ be the (non-singular) Hilbert scheme of ideal sheaves ${\mathcal I}_Z$ with $l({\mathcal O}_Z) = n$. 

\medskip

Our goal here is to get a useful lower bound for $m = \tan(\theta)$ so that the sector: 
\[ \left( I({\mathcal I}_Z) + z \ | \ \arg(z) < \theta\right))  \cap (x > \frac 12 y^2) \] 
consists of stability conditions $(x,y)$ for which {\bf all} the ideal sheaves are $(x,y)$-stable. 
This will translate in \S 7 to a good estimate for the size of the slice of the ample cone of $S^{[n]}$ in the space generated by 
 $(K_S + D)^{[n]}$ and $E$ i.e. Theorem 7.2.

\medskip

\nt Note. The line passing through $(-n,0) = I({\mathcal I}_Z)$ and tangent to the parabola
\[ \ \mbox{(i) \ has slope} \ \frac 1{\sqrt {2n}} \ \mbox{and (ii) \ meets the parabola at} \ \left( n, \sqrt{2n} \right)  \  \] 
so each wall has slope $< \frac 1{\sqrt{2n}}$ and intersects the parabola at $(x_0,y_0)$ with $x_0 < n$. 
This is an upper bound for the slope of a wall. The lower bound is more subtle. 

\medskip

\nt {\bf Preliminaries.} A {\it witness} for a wall  $W$ through $I({\mathcal I}_Z) = (-n,0)$ is a sequence: 
\[ 0 \rightarrow F \rightarrow {\mathcal I}_Z \rightarrow E \rightarrow 0 \] 
of semi-stable objects of the same $(x,y)$-slope, so that $W$ passes through $I(F)$ and, together with $\arg(z) = 0$ bounds a sector in which 
all ideal sheaves ${\mathcal I}_Z$ are stable. 
We conclude as before that: 

\medskip

$\bullet$ $F$ is a torsion-free coherent sheaf with $c_1(F) \cdot H < 0$. 

\medskip

$\bullet$ $I(F)$ is between $(-n,0)$ and the parabola.   

\medskip

$\bullet$ if $r = \rk(F) > 1$, then $I(E)$ is on the opposite side of the parabola and
\[ -\mu(E) = \frac{c_1(E)\cdot H}{r-1} \ \mbox{and} \ -\mu(F) = \frac{c_1(E)\cdot H}{r} \] 
so that in this case the {\it ratio} of the $y$-coordinates of $I(E)$ and $I(F)$ is $r/(r-1)$. 

\medskip

\nt {\bf Proposition 5.1.} If the witness $F \subset {\mathcal I}_Z$ for a wall $W$ has rank $r > 1$, then
\[ \frac 2{3\sqrt n} \ \mbox{
is a lower bound for the slope of $W$} \]

{\bf Proof.} Let $y_1,y_2$ be the $y$-coordinates of $W$ intersect the parabola. Then 
\[ -\mu(E) \le y_1 < y_2 \le -\mu(F) \] 
and so their ratio satisfies 
${y_2}/{y_1} \le{r}/{(r-1)}$.
On the other hand, from: 
\[ y^2 - \frac 2m y + 2n = 0 \] 
we conclude that
$y_1y_2 = 2n$ and $y_1 + y_2 = \frac 2m$
which gives us:
\[ \frac 2m \le \sqrt{2n} \cdot \left( \sqrt{ \frac r{r-1}} + \sqrt{\frac{r-1}r }\right) \]  
and a lower bound for $m$. 
When $r = 2$, this is the lowest lower bound, namely  
\[ m \ge \frac 23 \cdot \frac 1{\sqrt n}  \]

\nt {\bf Example.} When $n = 1$ and $S = \PP^2$, the wall for ${\mathcal I}_p$  is witnessed by: 
\[ 0 \rightarrow {\mathcal O}_{\PP^2}(-1)^2 \rightarrow {\mathcal I}_p \rightarrow {\mathcal O}_{\PP^2}(-2)[1] \rightarrow 0 \]  
and indeed, this realizes the lower bound for higher rank walls, as 
\[ I(F) = \left( \frac 12, 1\right) \ \mbox{and} \ I(E) = \left(2,2\right)\ \mbox{and} \ m = \frac 23 \] 

Next, we turn to the rank one walls for $S^{[n]}$. In general, a wall is witnessed by: 
\[ {\mathcal I}_W(-C) \subset {\mathcal I}_Z \] 
for a nonzero effective curve $C$ through a subset of the points of $Z$. 

\medskip

But if there is such a wall  then there is a wall of lower slope witnessed by: 
\[ {\mathcal O}_S(-C) \subset {\mathcal I}_{Z'} \] 
for $Z' \subset C$. So although witnesses of the first type do
exist, for the purpose of finding the {\it first} wall  it suffices to consider only the walls of the latter type. 

\medskip

\nt {\bf Proposition 5.2.} Suppose $H = D/||D|$ for an integral divisor $D$ such that: 
\[ D \cdot C > C^2 + 2n \ \mbox{for all} \ C \ \mbox{satisfying} \ -2n < C^2 < 2n \] 
Then the slopes of all rank one walls are bounded below by $2/||D||$. 

\medskip

{\bf Proof.} The condition for a wall with witness ${\mathcal O}_S(-C)$ is that: 
\[ I({\mathcal O}_S(-C)) = \left( \frac{C^2}2, C \cdot H \right) \ \mbox{be between $(-n,0)$ and the parabola}  \] 

The slope of the wall is then: 
\[ \frac{ C\cdot H}{\frac{C^2}2 + n} = \frac{2 C \cdot H} {C^2 + 2n} = \frac {2C \cdot D}{||D|| (C^2 + 2n)}> \frac 2{||D||} \] 
when $C\cdot D > C^2 + 2n$. This only needs to be checked 
for $-2n < C^2 < 2n$ since $(n,\sqrt{2n})$ is the point of tangency of the line through $(-n,0)$ and the parabola.

\medskip

Putting the two Propositions together, we obtain: 

\medskip

\nt {\bf Corollary 5.3.} When $D$ is an ample divisor on $S$ satisfying:
\[ D \cdot D > 9n \ \mbox{and} \ C \cdot D > C^2 + 2n \ \mbox{for all} \ C^2 < n \] 
then each ${\mathcal I}_Z$ is stable up to and including the wall through $I({\mathcal I}_Z)$ with slope $2/||D||$.  

\medskip

\nt {\bf Proof.} Except for the fact that it is sufficient to check $C^2 < n$ instead of $C^2 < 2n$, this is an immediate consequence of the two 
Propositions. 

\medskip

But let
$C^2 = \alpha n \ \mbox{for} \ 0 < \alpha < 2$ and suppose
\[  D \cdot C \le C^2 + 2n  = (\alpha + 2)n = \left( 1 + \frac 2\alpha\right)C^2 \] 

Then it follows from the Hodge Index Theorem that when $1 \le \alpha < 2$, 
\[ D^2 \le \left( 1 + \frac 2\alpha\right) D\cdot C \le \left( 1 + \frac 2\alpha\right)^2 C^2 = \frac{(\alpha + 2)^2}{\alpha} n \le 9n  \] 
 So $D^2 > 9n$  implies we need only check $D \cdot C > C^2 + 2n$ for $\alpha < 1$, i.e. $C^2 < n$.  \qed

\medskip

We will need in \S 6 to know when each ${\mathcal I}_p$ is stable along a hypothetical wall:
\[ W(a) = \left( y = \frac 2{||D||}x - \frac 2{||D||} a\right) \cap \left( x > \frac 12y^2\right)  \ \mbox{for} \ -2 \le a < -1 \] 
for the Drinfeld family.
Again, we can separate into higher rank and rank one. 

\medskip

\nt (Higher Rank) It suffices to know that the hypothetical wall satisfies:
\[ \frac 2{||D||} \cdot \frac 12 - \frac 2{||D||}a \le 1 \]  
to know it passes  below the open line segment joining $(\frac 12, 1)$ and $(2,2)$. This gives: 
\[ D^2 > (1 - 2a)^2 \] 
and in particular, 
$D^2 > 16 \ \mbox{when} \ a = -\frac 32 \ \mbox{and} \ D^2 > 25 \ \mbox{when} \ a = -2$.

\medskip 

\nt (Rank One) Given that the hypothetical wall lies below $(\frac 12,1)$ with the inequalities above, it suffices to 
test curves with $C^2 = -1, 0$, and to check that the line through 
$(-1,0)$ and $I({\mathcal O}_S(-C))$ passes through or above $(\frac 12, 1)$, i.e. that 
\[  m = \frac{C\cdot H}{\frac 12 C^2 + 1} \ge -\frac{2a}{||D||}  \] 
when $I({\mathcal O}_S(-C))$ is between $(-1,0)$ and the parabola. That is, 
\[ C \cdot D \ge -a(C^2 + 2) \ \mbox{when} \ C^2 = -1,0 \]  

Finally, the argument for ${\mathcal I}^\vee_p(-D)[1]$ is completely analogous, and we have: 

\medskip

\nt {\bf Corollary 5.4} (a) The stability sectors for ${\mathcal I}_p$ and ${\mathcal I}_p^\vee(-D)[1]$ contain $W(-\frac 32)$ when
 \[ D^2 > 16 \ \mbox{and} \  C \cdot D \ge \frac 32(C^2 + 2)  \ \mbox{for} \ C^2 = -1,0 \] 

(b) The stability sectors for ${\mathcal I}_p$ and ${\mathcal I}_p^\vee(-D)[1]$ contain $W(-2)$ when
\[ D^2 > 25 \ \mbox{and} \  C \cdot D \ge 2(C^2 + 2)  \ \mbox{for} \ C^2 = -1,0 \] 

\nt Note. The  $C^2 = 1$ case is ruled out by the inequality on $D^2$. 

\section{Stability Conditions and the Drinfeld Family} 

There is a  ``relative Serre'' family $C_{|K_S + D - p - q|^\vee}$  of extensions: 
\[ 0 \rightarrow {\mathcal I}_q^\vee(-D)[1] \rightarrow C_\epsilon \rightarrow {\mathcal I}_p \rightarrow 0 \] 
fibered over the product $S \times S$ with fibers:
\[ \PP(\mbox{H}^0(S, \omega_S(D) \otimes {\mathcal I}_p \otimes {\mathcal I}_q)) = |K_S + D - p - q|^\vee \]  

This intersects the Drinfeld family along the exceptional divisor $E$, which embeds in the 
relative Serre family via 
\[ \PP(\mbox{H}^0(S, \omega_S(D) \otimes {\mathcal I}_p^2)) \subset |K_S + D - p - p|^\vee \]

Thus, the objects of the  Drinfeld family are all $(x,y)$-(semi)stable outside $W_1$ if:  

\medskip

(i) $D^2 > 10$ and $D \cdot C > C^2 + 2 \ \mbox{for} \ C^2 = -1,0$.

\medskip

(ii) All objects  $C_\epsilon$ of the original Serre family for $\epsilon \not\in S$  are $(x,y)$-(semi)stable. 

\medskip

(iii) All objects $C_\epsilon$ of the relative Serre family are $(x,y)$-(semi)stable. 

\medskip

So we assume (i) and go to work on (ii) and (iii). 

\medskip

\nt {\bf Original Serre.}  If $C_\epsilon$ is produced by the Serre family with $\epsilon \not\in S$, then: 

\medskip

(Higher Rank) We've already made this computation in \S 4. If 
\[ D^2 > (-8a)\left( \frac{r^2}{r^2 -1}\right)  \] 
then $W(a)$ is not witnessed by a rank $r$ (or more) sheaf, and in particular, if 
\[ D^2 > (-a) \frac{32}3 \] 
then $W(a)$ is not witnessed in rank $> 1$. 

\medskip

(Rank One) The condition not to be witnessed by ${\mathcal I}_Z(-C) \subset C_\epsilon$ is also done:  
\[ D\cdot C > C^2  - 2a \ \mbox{jumps when $a = -\frac 32$ and when $a = -2$}\]

\nt {\bf Relative Serre} The genesis wall is $W_1$. 
If  $E \subset C_\epsilon$ witnesses  the next wall, with: 
\[ 0 \rightarrow {\mathcal I}_p^\vee(-D)[1] \rightarrow C_\epsilon \rightarrow {\mathcal I}_q \rightarrow 0 \]  
we would first like to conclude, as in \S 4, that 
$E \cap {\mathcal I}_p(-D)^\vee[1] = 0$ 
 and therefore that 
$E \subset {\mathcal I}_q$ is a torsion-free coherent sheaf. 
From the cohomology sequence:
\[ 0 \rightarrow {\mathcal O}_S(-D)[1] \rightarrow {\mathcal I}_p^\vee(-D) \rightarrow \CC_p \rightarrow 0 \] 
we see that
$ {\mathcal H}^{-1}(E) \subset {\mathcal O}_S(-D)$ and then as before that either
\[ {\mathcal H}^{-1}(E) = 0 \ \mbox{(and $E$ is a torsion-free sheaf), or} \ {\mathcal H}^{-1}(E) = {\mathcal O}_S(-D) \] 

In the latter case, ${\mathcal H}^0(E) \subset {\mathcal H}^0(C_\epsilon)$ in ${\mathcal A}_{-y}$, and from the sequence:
\[ 0 \rightarrow \CC_p \rightarrow {\mathcal H}^0(C_\epsilon) \rightarrow {\mathcal I}_q \rightarrow 0 \] 
it follows that either:

\medskip

(i) $\CC_p \subset {\mathcal H}^0(E)$ and ${\mathcal I}_p^\vee(-D)[1] \subset E$ or else

\medskip

(ii) $\CC_p \cap {\mathcal H}^0(E) = 0$, in which case replacing $E$ with $E' = 
\ker(C_\epsilon \rightarrow {\mathcal H}^0(C_\epsilon)/{\mathcal H}^0(E))$
results in a subobject of larger $(x,y)$-slope, and a contradiction.  

\medskip

We may draw a contradiction as in \S 4 to the alternative $I_Z^\vee(-D)[1] \subset E$ when we use the  additional 
inequalities of Corollary 5.4 to conclude that ${\mathcal I}_p$ is stable along the wall. Thus with these inequalities,
we may conclude that an $E$ withnessing the wall is a torsion-free sheaf. At that point, the 
argument proceeds exactly as with the original Serre family, and we obtain: 

\medskip

\nt {\bf Proposition 6.1.} (a) if $D^2 > 16$ and $D \cdot C > C^2 + 3$ when $C^2 = -1,0$ then the stability region for the Drinfeld family contains the wall
$W(-\frac 32)$. 

\medskip

(b) if $D^2 > 25$ and $D \cdot C > C^2 + 4$, the stability region contains the wall $W(-2)$. 

\medskip

\section{The Determinant Construction} 

The construction of Bayer-Macr\`i \cite{BM3} associates a linear map: 
\[ \gamma_{(x,y)}: N_1(X) \rightarrow \RR  \] 
to each family $F_X$ of objects of ${\mathcal A}_{-y}$ over a variety $X$. It is defined by: 
\[ \gamma_{(x,y)}(C) := - \mathfrak {Im} \left( \frac { Z_{(x,y)}(R\pi_* F_C)}{Z_{(x,y)}(F_p)} \right) \] 
for $p\in C \subset X$ and satisfies the following properties: 

\medskip

$\bullet$  $\gamma_{(x,y)}(C) \ge 0$ if $F_C$ is a family of $(x,y)$-semistable objects, and among these

\medskip

$\bullet$ $\gamma_{(x,y)}(C) = 0$ only when the $(x,y)$-associated gradeds $\mbox{Ass}(F_p)$ do not vary. 

\medskip

$\bullet$ the ray spanned by $\gamma_{(x,y)}$ is constant along walls for the invariants $\ch(F_p)$. 

\medskip

We now follow a strategy outlined by Bayer and Macr\`i to prove the Theorems. 

\medskip

\nt {\bf Theorem 7.1.} The nef cone of the blow-up $Y$ of $|K_S + D|^\vee$ along $S$ 

\medskip

(i) contains $3H - E$ when $D^2 \ge 16$ and $D \cdot C \ge C^2 + 3$ for all $C^2 = -1,0$. 

\medskip

(ii) contains $3H - E$ in the {\it interior} when the inequalities are strict. 

\medskip

(iii) contains $2H - E$ when $D^2 \ge 25$ and $D \cdot C \ge C^2 + 4$ for $C^2 = -1,0$. 

\medskip

(iv) contracts the variety of secant lines onto $S^{[2]}$ when the inequalities are strict. 

\medskip

{\bf Proof.} We compute the determinant divisor class for the Drinfeld family: 

\medskip

(a) Let $l \subset |K_S + D|^\vee$ be a line not meeting $S$ (we assume codim$_{|K_S + D|^\vee}S \ge 2$).
Then via the restriction $C_l$ of the Serre family, we obtain: 
\[ \gamma_{(x,y)}(l) = - \mathfrak{Im} \frac{ Z_{(x,y)}(-1,2D,-D^2)}{Z_{(x,y)}(0,D,-\frac {D^2}2)} = (-a) \left( \frac{4||D||}{D^2(D^2 + 4)} \right)  \]  
for $(x,y) \in W(a)$, the wall through $(a,0)$ with slope $2/||D||$.   

\medskip

(b) Let $l_p$ be a line $l_p$ in a fiber of $\pi:E \rightarrow S$. Then: 
\[ \gamma_{(x,y)}(l_p) = - \mathfrak{Im} \frac{ Z_{(x,y)}(-1,2D,-D^2 + 1)}{Z_{(x,y)}(0,D,-\frac {D^2}2)} = -(a + 1) \left( \frac{4||D||}{D^2(D^2 + 4)} \right) \] 
from which we conclude that
$\gamma_{(x,y)} \ \mbox{is proportional to} (-a)H + (-a + 1)E$ on $W(a)$, i.e. to
$3H - E \ \mbox{along} \ W(-3/2) \  \mbox{and} \  2H - E \ \mbox{along} \ W(-2)$.

\medskip

Together with Proposition 6.1, this gives (i)-(iii). 

\medskip

For (iv), we consider the strictly semi-stable 
objects of $Y$ along the wall $W(-2)$. This is analogous to the locus $S \subset |K_S + D|^\vee$ of strictly semi-stable 
objects of the Serre family along the wall $W(-1)$. By Corollary 5.3 and 5.4, we may conclude that all of the objects ${\mathcal I}_p, {\mathcal I}_p^\vee(-D)[1]$ as well as 
${\mathcal I}_Z, {\mathcal I}_Z^\vee(-D)[1]$ for $l(Z) = 2$ are stable along the wall, and then liftability of the ideal sheaf of $Z \in S^{[2]}$: 
\[ \begin{array}{cccccccccc} 
 {\mathcal I}_Z \\  \downarrow & \searrow \\ C_\epsilon & \rightarrow & {\mathcal O}_S & \rightarrow & {\mathcal O}_S(-D)[1] & \rightarrow \\
 \end{array} \] 
flags the elements $C_\epsilon$ (for $\epsilon \not \in S$) that are strictly semi-stable. In this case, however, the object $C_\epsilon$ is not uniquely 
determined by the liftability. Instead, all points of the secant (or tangent) line spanned by $Z$ (other than points of $S$) admit the lift. 

\medskip

This is the non-singular complement $\Sigma_1(S) - S$ of the variety of 
secant/tangent lines to $S$, whose closure in $Y$ gives the strictly semi-stable locus fibered over the moduli space $S^{[2]}$ 
of $(x,y)$-stable objects of class $(1,0,-2)$ (for $(x,y) \in W(-2)$).

\medskip

Next, we turn our attention to the Hilbert schemes more generally. 
Let ${\mathcal I}_{\mathcal Z}$ be the family of ideal sheaves over $S^{[n]}$. In this case the walls are the line segments (in the stability region) through
$(-n,0)$ of slope $m > 0$. We now compute the determinant divisor class (up to a scaling factor) along the wall $W(m)$ defined by: 
\[ y = m(x + n) \] 

(a) Fix $Z_0$ a reduced divisor of length $n-1$ and a curve $C \subset S - Z_0$ and consider the family of reduced divisors $F_C$ indexed by $C$ given by: 
$F_p = \{ p \cup Z_0 \}$.  
Then the family is an ideal in $C \times S$ satisfying: 
\[ 0 \rightarrow F_C \rightarrow {\mathcal I}_{Z_0 \times S} \rightarrow {\mathcal O}_\Delta \rightarrow 0 \] 
where $\Delta \subset C \times C$ is the diagonal. The derived push-forward of the family satisfies: 
\[ \ch(R\pi_* F_C) = (1 - g)\ch({\mathcal I}_{Z_0}) - \ch({\mathcal O}_C) = (1-g)(1, 0, -n) + (0,-C, \frac {C^2}2 + (1-g)) \]
while of course
$\ch(F_p) = (1,0,-n)$.  
We may simplify using the adjunction formula:
\[ (1 - g) = - \frac{C \cdot (C + K_S)}2 \] 
to obtain: 
\[ - \mathfrak{Im} \frac{Z_{(x,y)}(R\pi_* F_C)}{Z_{(x,y)}(F_p)} = \frac {y \cdot \frac{C \cdot K_S}2 + (C\cdot H)(n+x)}{(n + x)^2 + y^2}\] 
and then substituting $y = m(x + n)$ we further obtain:
\[ \gamma_{(x,y)}(C) = \left(\frac{m}{2(n+x)(1 + m^2)}\right)  \left( C \cdot (K_S + \frac {2D}{m||D||})\right)  \] 

(b) Next, let $Z_1$ contain one point $p\in C$ and $F_{C_p}$ the resulting family, fitting in: 
\[ 0 \rightarrow F_{C_p} \rightarrow {\mathcal I}_{Z_1 \times S} \rightarrow {\mathcal O}_\Delta(-p) \rightarrow 0 \] 
This yields: 
\[ \gamma_{(x,y)}(C_p) = \left(\frac{m}{2(n+x)(1 + m^2)}\right) \left( C \cdot (K_S + \frac{2D}{m||D||})   - 1 \right) \] 
and since $C_p \cdot E = 1$ because of the single point of the family with a non-reduced scheme, it follows that the determinant divisor class is proportional to
$( K_S + D) ^{[n]} - E$ 
along the line of slope $m = 2/||D||$. With Corollary 5.3, this gives Theorem 7.2.

\medskip

\nt Mathematics Department, University of Utah

\medskip

\nt Aaron Bertram, bertram@math.utah.edu

\nt Jonathon Fleck, fleckjem@math.utah.edu

\nt Liebo Pan, lpan@math.utah.edu

\nt Joseph Sullivan, jsullivan@math.utah.edu

\end{document}